\documentclass[11pt]{article}
\usepackage{amsmath, amsfonts,amssymb}
\usepackage{amsthm}
\theoremstyle{plain}
\newtheorem{Theorem}{Theorem}
\title{The circular $k$-partite crossing number of $K_{m,n}$}
\author{Adrian Riskin\\
Department of Mathematics\\
Mary Baldwin College\\
Staunton, VA  24401\\
ariskin@mbc.edu}
\begin{document}
\maketitle
\begin{abstract}
We define a new kind of crossing number which generalizes both the bipartite crossing number and the outerplanar 
crossing number.  We calculate exact values of this crossing number for many complete bipartite graphs and also
give a lower bound.
\end{abstract}

\section{Preliminaries}

The \textit{bipartite crossing number} of a bipartite graph $G$ was defined by Watkins in [W] to be the
minimum number of crossings over all bipartite drawings of $G$.  A \textit{bipartite drawing} of 
bipartite $G$ is one in which the vertices of the parts $V_{1}$ and $V_{2}$ are placed respectively on two
distinct parallel lines and then the edges of $G$ are drawn as straight line segments joining appropriate pairs
of vertices.  The calculation and estimation of this number are of interest to those who study VLSI design, 
graph drawing algorithms, and/or topological graph theory.  See [N] for a bibliography on the topic as well
as some of the few known exact results.

The \textit{outerplanar crossing number} of a graph $G$, also known as the \textit{circular} or
\textit{convex crossing number} of $G$, was defined by Kainen in [K] to be the minimum number of 
crossings taken over all plane drawings of $G$ where the vertices lie on a circle and the edges are chords
of that circle.  The calculation and estimation of this number are of interest to the same audience as the bipartite
crossing number.  See [C] and [F] for an introduction and bibliography, and [F] and [R] for the few known
exact results.

In this paper we introduce a notion, the \textit{circular k-partite crossing number} of a $k$-partite
graph $G$, which generalizes both of these definitions.  A \textit{circular k-partite drawing} of
$k$-partite $G$ is constructed as follows: partition a circle into $k$ segments of arc.  Place the vertices
of the $i^{th}$ part into the $i^{th}$ segment of arc and then add the edges as chords of the circle with
the usual proviso that no more than two edges should meet at a crossing.  The \textit{circular k-partite
crossing number} of $G$, denoted by $cpr_{k}(G)$, is the minimum number of crossings taken over all
circular $k$-partite drawings, all possible assignments of vertices to parts, and all numberings of the parts of $G$.
Note that if $G$ has $p$ vertices then $G$ is $p$-partite, and $cpr_{p}(G)=\nu_{1}(G)$, the familiar
outerplanar crossing number.  Furthermore, if $G$ is bipartite then $cpr_{2}(G)=bcr(G)$, the familiar
bipartite crossing number.

Note that we prepend the modifier ``circular'' to our epithet for $cpr_{k}(G)$ to distinguish it from the 
quite different $k$-partite crossing number, also known as the $k$-layer crossing number.  Finally, we
will have occasion to mention the following:

\begin{Theorem}\label{outerplanar}
If $m|n$ then the outerplanar crossing number of $K_{m,n}$ is
\begin{equation*}
\frac{1}{12}n(m-1)(2mn-3m-n)
\end{equation*}
\end{Theorem}

\noindent A proof of this may be found in [R].

\section{Results}

Note that $K_{m,n}$ is $k$-partite for $2 \leq k \leq m+n$.  Our goal is to determine $cpr_{k}(K_{m,n})$
for each $k$ in this range.  To that end we observe that if $m \leq n$ then $cpr_{k}(K_{m,n})=cpr_{2m}
(K_{m,n})$ for all $2m \leq k \leq m+n$, and that $cpr_{2k}(K_{m,n})=cpr_{2k+1}(K_{m,n})$ since
every circular $(2k+1)$-partite drawing of $K_{m,n}$ is also a $2k$-partite drawing and vice-versa.  Our
result, which generalizes Theorem \ref{outerplanar}, is the following:

\begin{Theorem}\label{main}
\begin{equation*}
cpr_{2k}(K_{m,n}) \geq {m \choose 2}{n \choose 2}-\frac{(k^{4}-k^{2})m^{2}n^{2}}{12k^{4}}
\end{equation*}
with equality when $k|m$ and $k|n$.
\end{Theorem}

\begin{proof} Let $D$ be a circular $2k$-partite drawing of $K_{m,n}$.  We denote the number of crossings in
$D$ by $cpr_{2k}(D)$.  Let $M$ and $N$ be the two parts of $K_{m,n}$, with $|M|=m$ and $|N|=n$.  We refer to
the vertices in $M$ as pink and to those in $N$ as black.  Let the $2k$ segments of arc be labeled consecutively 
$M_{1}, N_{1}, M_{2}, N_{2}, \dots , M_{k}, N_{k}$, and let $|M_{i}|=m_{i}$ and $N_{i}=n_{i}$.

Now suppose that the $2k$-partite sets fail to alternate colors.  In this case $D$ is a circular $2j$-partite drawing of
$K_{m,n}$ for some $j < k$ where the $2j$-partite sets do alternate colors.  Since
\begin{equation*}
{m \choose 2}{n \choose 2}-\frac{(k^{4}-k^{2})m^{2}n^{2}}{12k^{4}}
\end{equation*}
is a strictly decreasing function of $k$ for $k\geq 1$, if we prove the theorem for drawings where the partite sets do 
alternate colors we will have proved it also for drawings where they do not.

Let $u_{1}$ and $u_{2}$ be distinct pink vertices and let $v_{1}$ and $v_{2}$ be distinct black vertices.  Let
$u_{j}\in M_{i_{j}}$ and $v_{j} \in N_{i_{j}}$.  Then the vertices $u_{1}, u_{2}, v_{1}$, and $v_{2}$ determine
a crossing unless $M_{i_{1}}$ and $M_{i_{2}}$ separate $N_{i_{1}}$ and $N_{i_{2}}$ on the boundary of 
the circle.  If $1 \leq i < j \leq m$, then $M_{i}$ and $M_{j}$ separate $(j-i)(k-(j-i))$ distinct pairs of black partite 
sets from one another.  Explicitly, $N_{i}, \dots, N_{j-1}$ are separated from $N_{j}, \dots, N_{i-1}$ where the subscripts
are, naturally, read modulo $k$.  Suppose now that $u_{1} \in M_{i}$ and $u_{2} \in M_{j}$ with $1 \leq i < j \leq k$.
Then there are 
\begin{equation*}
\sum_{\substack{i \leq s \leq j-1\\
                           j \leq t \leq i-1}}
               m_{i}m_{j}n_{s}n_{t}            
\end{equation*}
choices of $v_{1}$ and $v_{2}$ which do not determine crossings.  Therefore
\begin{equation}\label{cprD}
cpr_{2k}(D)={m \choose 2}{n \choose 2}-
\sum_{\substack{1\leq i \leq k-1\\
                                i+1 \leq j \leq k}} \quad
                    \sum_{\substack{i \leq s \leq j-1\\
                                                j \leq t \leq i-1}}
         m_{i}m_{j}n_{s}n_{t} 
  \end{equation}
  Clearly $m_{i}m_{j}n_{s}n_{t} \leq \frac{m^{2}n^{2}}{k^{4}}$, and the lower bound follows from this (well, this
  and a considerable amount of algebra).  When $k|m$ and $k|n$, equality is obtained by distributing the vertices so that 
  there are $\frac{m}{k}$ in each $M_{i}$ and $\frac{n}{k}$ in each $N_{i}$.
\end{proof}

\noindent Note that we can obtain Theorem \ref{outerplanar} from Theorem \ref{main} by substituting $k=m$ in 
the case where $m|n$.  Note also that the double sum in \eqref{cprD} is maximized when the values of $m_{i}$ and
$n_{i}$ are as evenly distributed as possible, that is, when $m-k\left\lfloor\frac{m}{k}\right\rfloor$ of the $m_{i}$'s
are equal to $\left\lceil\frac{m}{k}\right\rceil$ and the other $k-m+
k\left\lfloor\frac{m}{k}\right\rfloor$ of them are equal to $\left\lfloor\frac{m}{k}\right\rfloor$
and likewise for the $n_{i}$'s.  It is possible to use \eqref{cprD} to obtain an exact expression for the value of 
$cpr_{2k}(K_{m,n})$ even when $k\nmid m$ or $k \nmid n$.  Some care must be taken to arrange the different values
as evenly as possible for values of $k\geq 3$. As a simple example, for $k=2$ the exact result is:
\begin{equation*}
cpr_{4}={m \choose 2}{n \choose 2} - \left\lceil\frac{m}{2}\right\rceil\left\lfloor\frac{m}{2}
\right\rfloor\left\lceil\frac{n}{2}\right\rceil\left\lfloor\frac{n}{2}
\right\rfloor
\end{equation*}
\bigskip

\section*{Bibliography}
\begin{description}

\item[[C]] Czabarka, E; S\'ykora, O.; Sz\'ekely, L.A.; Vrt'o, I. Outerplanar crossing numbers, the circular arrangement
problem, and isoperimetric functions.  Electronic J. Combinatorics 11(2004).

\item[[F]] Fulek, R.; He, H.; S\'ykora, O.; Vrt'o, I. Outerplanar crossing numbers of 3-row meshes, Halin graphs, 
and complete $p$-partite graphs.  Lecture Notes in Computer Science. 3381(2005) 376-379

\item[[K]] Kainen, P.C. The book thickness of a graph II.  Cong. Numer. 71(1990)

\item[[N]] Newton, M.C.; S\'ykora, O.; U\v zovi\v c, M.; Vrt'o, I.  New exact results and bounds for bipartite 
crossing numbers of meshes.  Lecture Notes in Computer Science. 3383(2005) 360-370

\item[[R]] Riskin, A. On the outerplanar crossing numbers of $K_{m,n}$.  Bull. Inst. Comb. Appl. 39(2003) 16-20.

\item[[S]] Shahrokhi, F.; Sz\'ekely, L.A.; S\'ykora, O.; Vrt'o, I. The book crossing number of a graph. J. Graph
Theory 21(1996) 413-424.

\item[[W]] Watkins, M.E.  A special crossing number for bipartite graphs: a research problem.  Ann. N.Y. Acad. 
Sci. 175(1970) 405-410.

\end{description}
\end{document}